\documentstyle[12pt,amsfonts,amssymb,leqno]{article}
\newfont\got{eufm10}

\newtheorem{proposition}{Proposition}[section]
\newtheorem{thm}[proposition]{Theorem}

\newtheorem{lemma}[proposition]{Lemma}

\renewcommand{\thefootnote}{\alph{footnote}}

\newcounter{secnum}

\begin{document}
\setcounter{section}{+0}

\begin{center}
{\Large \bf Cardinal arithmetic and Woodin cardinals}
\end{center}
\begin{center}
{\large Ralf Schindler}
\renewcommand{\thefootnote}{arabic{footnote}}
\end{center}
\begin{center} 
{\footnotesize
{\it Institut f\"ur Formale Logik, Universit\"at Wien, 1090 Wien, Austria}} 
\end{center}

\begin{center}
{\tt rds@logic.univie.ac.at}

{\tt http://www.logic.univie.ac.at/${}^\sim$rds/}\\
\end{center}


\begin{abstract}
\noindent Suppose that there is a measurable cardinal.
If $2^{\aleph_1} < \aleph_\omega$, but $\aleph_\omega^{\aleph_0} >
\aleph_{\omega_1}$, then there is an inner model
with a Woodin cardinal. This essentially answers a question of Gitik and Mitchell
(cf.~\cite[Question 5, p.~315]{moti-bill}).
\end{abstract}

We refer the reader to \cite{uri-menachem} for an introduction to cardinal
arithmetic and to Shelah's pcf theory (cf.~also \cite{burke-magidor}).
The perhaps most striking result of Shelah's in cardinal arithmetic is that if
$\aleph_\omega$ is a strong limit cardinal then $$2^{{\aleph_\omega}} < {\rm min}(
\aleph_{(2^{\aleph_0})^+} , \aleph_{\omega_4}).$$
Magidor was the first one to produce a model of set theory in which
the ${\sf GCH}$ holds below $\aleph_\omega$, but $2^{\aleph_\omega} =
\aleph_{\omega+2}$ (cf.~\cite{menachem1}, \cite{menachem2}). It is now known how to
produce models in which there are arbitrarily large countable gaps between
$\aleph_\omega$ and $2^{\aleph_\omega}$, while 
the ${\sf GCH}$ holds below $\aleph_\omega$ (cf.~for instance \cite{moti-menachem}).
A strong cardinal is more than enough for this purpose.
In fact, many equiconsistencies are known. 
We refer the reader to \cite{moti} and \cite{bill} and to the references given there. 

It is, however, open if it is possible to have
that $\aleph_\omega$ is a strong limit cardinal, but $2^{\aleph_\omega} >
\aleph_{\omega_1}$ (cf.~\cite[Section 7, Problem 1]{moti}).
This problem is just one of the key open problems of pcf theory in disguise.
Gitik and Mitchell have shown that a strong cardinal is not enough for 
producing such a model:

\begin{thm}\label{moti-bill-thm} {\em (\cite[Theorem 5.1]{moti-bill})}
If $2^{\aleph_0} < \aleph_\omega$ and $\aleph_\omega^{\aleph_0} >
\aleph_{\omega_1}$ then
there is a sharp for a model with a strong cardinal.
\end{thm}

The purpose of this note is to prove the following theorem which in a certain sense
improves Theorem \ref{moti-bill-thm}. It will say that you will need at least a Woodin
cardinal in order to produce a model in which $\aleph_\omega$ is a strong limit 
cardinal, but $2^{\aleph_\omega} >
\aleph_{\omega_1}$. As far as we know this is the first statement in cardinal
arithmetic which is known to practically imply the consistency of a Woodin cardinal.

\begin{thm}\label{main-thm0}
If $2^{\aleph_1} < \aleph_\omega$, $\aleph_\omega^{\aleph_0} >
\aleph_{\omega_1}$,
and there is a measurable cardinal then there is an inner model
with a Woodin cardinal.
\end{thm}

Our proof of Theorem \ref{main-thm0} will make use of Shelah's pcf theory. 
Specifically, we'll need the following theorems which are due to Shelah.
Recall that if ${\sf a}$ is a set of regular cardinals then
${\rm pcf}({\sf a})$ is the set of all possible cofinalities of
$\prod {\sf a} / U$ where $U$ is an ultrafilter on ${\sf a}$.

\begin{thm}\label{pcf1} {\em (\cite[Theorem 5.1]{burke-magidor})}
Let $2^{\aleph_0} < \aleph_\omega$. Then ${\rm pcf}(\{ \aleph_n \colon n<\omega \})
= \{ \kappa \leq \aleph_\omega^{\aleph_0} \colon \kappa$ is regular $ \}$.
\end{thm}

\begin{thm}\label{pcf2} {\em (\cite[Theorem 6.10]{burke-magidor})}
Let $2^{\aleph_0} < \aleph_\omega$. Let ${\sf d} \subset 
{\rm pcf}(\{ \aleph_n \colon n<\omega \})$ and $\mu \in {\rm pcf}({\sf d})$. 
There is then
some ${\sf d}' \subset {\sf d}$ such that ${\rm Card}({\sf d}') = \aleph_0$ and 
$\mu \in {\rm pcf}({\rm d}')$.
\end{thm}

We shall also need the following simple ``combinatorial'' fact, Lemma \ref{key-lemma}.
Let $\kappa$ and $\lambda$ be cardinals with $\lambda \leq \kappa$. 
$H_\kappa$ is the set of all sets which are hereditarily
smaller than $\kappa$, and $[H_\kappa]^\lambda$ is the set of all subsets of
$H_\kappa$ of size $\lambda$.
Recall that $S \subset [H_\kappa]^\lambda$ is called stationary if and only if
for all models
${\frak M} = (H_\kappa;...)$ of finite type and with universe $H_\kappa$ there is some
$X \in S$ such that $X$ is the universe of an elementary submodel of ${\frak M}$,
i.e., $(X;...) \prec {\frak M}$. We say that $S \subset [H_\kappa]^\lambda$ is
$*$-stationary if and only if
$S \cap \{ x \colon {}^\omega x \subset x \}$ is stationary, i.e.,
if for all models
${\frak M} = (H_\kappa;...)$ of finite type and with universe $H_\kappa$ there is some
$X \in S$ such that ${}^\omega X \subset X$ and
$X$ is the universe of an elementary submodel of ${\frak M}$. 
We let $NS_{\omega_1}$ denote the non-stationary ideal on $\omega_1$. 

\begin{lemma}\label{key-lemma}
Let $\kappa \geq 2^{\aleph_1}$ be regular, and let
$\Phi \colon [H_\kappa]^{2^{\aleph_1}} \rightarrow NS_{\omega_1}$.
There is then a pair $(C,S)$ such that $C$ is a closed unbounded subset of $\omega_1$,
$S$ is $*$-stationary in
$[H_\kappa]^{2^{\aleph_1}}$, and $C \cap \Phi(X) = \emptyset$ for all $X \in S$. 
\end{lemma}

{\sc Proof}. Suppose that for every club $C \subset \omega_1$ the set
$$\{ X \in [H_\kappa]^{2^{\aleph_1}} \colon C \cap \Phi(X) = \emptyset \}$$ is not
$*$-stationary. This means that for every club $C \subset \omega_1$ there is a model
${\frak M}_C$ of finite type and with universe $H_\kappa$ such that for every
$(X;...) \prec {\frak M}_C$ with ${}^\omega X \subset X$ and ${\rm Card}(X) =
2^{\aleph_1}$ we have that $C \cap \Phi(X) \not= \emptyset$.
As there are only $2^{\aleph_1}$ many subsets of $\omega_1$, 
there is a model ${\frak M}$ of finite type and with universe $H_\kappa$ such that
for every club $C \subset \omega_1$, if $(X;...) \prec {\frak M}$ is such that
$2^{\aleph_1} \subset X$ then $(X;...) \prec {\frak M}_C$ 
for every club $C \subset \omega_1$.
Pick $(X;...) \prec {\frak M}$ with ${}^\omega X \subset X$, ${\rm Card}(X) =
2^{\aleph_1}$, and $2^{\aleph_1} \subset X$.
We shall have that $C \cap \Phi(X) \not= \emptyset$ for every club $C \subset
\omega_1$, which means that $\Phi(X)$ is stationary. Contradiction! 
\hfill $\square$

\bigskip
Our proof of Theorem \ref{main-thm0} will use the core model theory of \cite{CMIP}.
The basic idea for its proof will be the following. We shall first use Theorems
\ref{pcf1} and \ref{pcf2} as well as Lemma \ref{key-lemma} for isolating
a ``nice'' countable set ${\sf d}' \subset \{ \aleph_{\alpha+1} \colon \alpha <
\omega_1 \}$ with $\aleph_{\omega_1+1} \in {\rm pcf}({\sf d}')$. We shall then use
a covering argument to prove that ${\rm pcf}({\sf d}') \subset \aleph_{\omega_1}$
yielding the desired contradiction. However, the covering argument plays the key role
in choosing the ``nice'' ${\sf d}'$ we start with.

\bigskip
{\sc Proof} of Theorem \ref{main-thm0}.
Suppose not. Let $\Omega$ be a measurable cardinal, and let $K$ denote
Steel's core model of height $\Omega$ (cf.~\cite{CMIP}). Let $\kappa >
\aleph_{\omega_1}$ be a regular cardinal.

Fix $X \prec H_\kappa$ for a while, where ${}^\omega X \subset X$.
Let $\pi \colon H_X \cong X \prec H_\kappa$ be such that $H = H_X$ is transitive.
Let ${\bar K} = {\bar K}_X = \pi^{-1}(K ||
\aleph_{\omega_1})$. We know by \cite{covering} that there is an $\omega$-maximal
normal iteration tree ${\cal T}$ on $K$ (of successor length) 
such that ${\cal M}_\infty^{\cal T} \trianglerighteq
{\bar K}$. Let ${\cal T} = {\cal T}_X$ be the shortest such tree.

If $E$ is an extender then we shall write $\nu(E)$ for the natural length of $E$
(cf.~\cite[p.~6]{FSIT}). We shall let $$\Phi(X) =
\{ \alpha < \omega_1 \colon \exists \beta+1 \in (0,\infty]_T \ ( {\rm
crit}(E_\beta^{\cal T}) < \aleph_\alpha^H \leq \nu(E_\beta^{\cal T}) ) \}.$$
We aim to apply Lemma \ref{key-lemma} to $\Phi$.

\bigskip
{\bf Claim 1.} $\Phi(X)$ is a non-stationary subset of $\omega_1$.
\bigskip

{\sc Proof.} Suppose that $S_0 = \Phi(X)$ is stationary. Let $S$ be the set of all
limit ordinals of $S_0$. $S$ is stationary, too. Let $F \colon S \rightarrow {\rm OR}$
be defined by letting $F(\alpha)$ be the least ${\bar \alpha} < \alpha$ such that
${\rm crit}(E_\beta^{\cal T}) < \aleph_{\bar \alpha}^H$, where 
$\beta+1 \in (0,\infty]_T$ is unique such that ${\rm
crit}(E_\beta^{\cal T}) < \aleph_\alpha^H \leq \nu(E_\beta^{\cal T})$. Let ${\bar S}
\subset S$ be stationary such that $F \upharpoonright {\bar S}$ is constant.
For $\alpha \in {\bar S}$ let $\beta(\alpha)$ be the unique $\beta$ such that ${\rm
crit}(E_\beta^{\cal T}) < \aleph_\alpha^H \leq \nu(E_\beta^{\cal T})$.
Using the initial segment condition \cite[Definition 1.0.4 (5)]{FSIT}
(cf.~also 
\cite[Definition 2.4]{deconstructing}) it is easy to see that we must have ${\rm
crit}(E_{\beta(\alpha)}^{\cal T}) = {\rm
crit}(E_{\beta(\alpha')}^{\cal T})$ whenever $\{ \alpha , \alpha' \} \subset {\bar S}$.
Hence $E_{\beta(\alpha)}^{\cal T} = E_{\beta(\alpha')}^{\cal T}$
whenever $\{ \alpha , \alpha' \} \subset {\bar S}$. Let us write $E$ for this unique
extender. We'll have to have $\nu(E) \geq \aleph_{\omega_1}^H$, so that $E$ cannot
have been used in ${\cal T}$. Contradiction! \hfill $\square$

\bigskip
We shall now define $\Phi \colon [H_\kappa]^{2^{\aleph_1}} \rightarrow NS_{\omega_1}$
as follows. Let $X \in [H_\kappa]^{2^{\aleph_1}}$. If $X \prec H_\kappa$ is such that 
${}^\omega X \subset X$ then we let $\Phi(X)$ be defined as above. Otherwise we set 
$\Phi(X) = \emptyset$. Let $(C,S)$ be as given by Lemma \ref{key-lemma}.\footnote{The
lift-up arguments which are to follow will be simplified by the assumption, which we
may make without loss of generality, that every element of $C$ is a limit ordinal.}
We let $${\sf d} = \{ \aleph_{\alpha+1} \colon \alpha \in C \}.$$
We know that $\aleph_{\omega_1+1} \in {\rm pcf}({\sf d})$, by \cite[Remark
1.8]{burke-magidor}. By Theorems \ref{pcf1} and \ref{pcf2} there is a countable ${\sf
d}'
\subset {\sf d}$ with $\aleph_{\omega_1+1} \in {\rm pcf}({\sf d}')$. We shall now 
derive a
contradiction by showing that ${\rm pcf}({\sf d}') \subset {\rm sup}({\sf d}')^++1$.

Let $\lambda = {\rm sup}({\sf d}')$. 

\bigskip
{\bf Main Claim.} For every $f \in \prod d'$ there is some $g \in K$, $g \colon
\lambda \rightarrow \lambda$ such that $g(\mu) > f(\mu)$ for all $\mu \in {\sf d}'$.

\bigskip
Suppose that the Main Claim holds. Then certainly $$\{ [g \upharpoonright {\sf d}']_U
\colon g \in K \wedge g \colon \lambda \rightarrow \lambda \}$$ is cofinal in $\prod
{\sf d}' / U$. But there are only $\leq \lambda^+$ many $g \in K$, 
$g \colon \lambda \rightarrow \lambda$, so that we must cartainly have
${\rm cf}(\prod {\sf d}' / U) \leq \lambda^+$. This contradiction proves Theorem
\ref{main-thm0}.

It therefore suffices to prove the Main Claim.
Fix $f \in \prod d'$ for the rest of this proof. 
By the choice of $(C,S)$ there is some $X \prec
H_\kappa$ with ${\rm Card}(X) = 2^{\aleph_1}$, ${}^\omega X \subset X$, 
$C \cap \Phi(X) = \emptyset$, and $f \in X$. 
Let $\pi \colon H = H_X \cong X \prec H_\kappa$, and let
${\cal T} = {\cal T}_X$ be as defined above. Then if
$\beta+1 \in (0,\infty]_{T}$ and $\aleph_{\alpha+1} \in {\sf d}'$ we do not have that
${\rm crit}(E_\beta^{\cal T}) < \aleph_{\alpha}^H \leq \nu(E_\beta^{\cal T})$.

Let $(D^{\cal T} \cap (0,\infty]_T) \cup \{ {\rm lh}({\cal T}) \}
= \{ \alpha_0+1 < ... < \alpha_N+1 \}$, where $0 \leq N < \omega$.
Let $\alpha_n^*$ be the $T$-predecessor of $\alpha_n+1$ for $0 \leq n \leq N$.
Let us pretend that $N > 0$ and $\alpha_0^* = 0$, i.e., that we 
immediately drop on $(0,\infty]_T$.\footnote{I.e., in this note we simply ignore the
possibility that we might have
$\alpha_0^* > 0$. In any event, it can be shown that if we were to
use Friedman-Jensen premice rather than Mitchell-Steel premice then 
the case that $\alpha_0^* > 0$ would not come up by an argument of \cite{habil}.}
Notice that ${\rm crit}(\pi) \leq \pi^{-1}((2^{\aleph_1})^+)$. Let us assume without
loss of generality that $\alpha \in C \Rightarrow \aleph_\alpha > (2^{\aleph_1})^+$. 

For each $\alpha \in C$ there is a least $\beta(\alpha) \in [0,\infty]_T$ 
such that
${\cal M}_{\beta(\alpha)}^{\cal T} || \aleph_\alpha^H = {\bar K} || \aleph_\alpha^H$.
Let $n(\alpha)$ be the unique $n<N$ such that
$\alpha_{n-1}^* < \beta(\alpha) \leq \alpha_n^*$.
Then $D^{\cal T} \cap (\beta(\alpha),\alpha_{n(\alpha)}^*]_T = \emptyset$ and 
$\rho_\omega({\cal M}^{\cal T}_{\alpha_{n(\alpha)}^*}) \leq \aleph_\alpha^H$.
Let $\eta(\alpha)$ be the least $\eta$ such that 
$\rho_\omega({\cal M}^{\cal T}_{\alpha_{n(\alpha)}^*} || \eta(\alpha)) 
\leq \aleph_\alpha^H$ 
(if $\beta(\alpha) < \alpha_n^*$ then $\eta(\alpha) = 
{\cal M}^{\cal T}_{\alpha_{n(\alpha)}^*} \cap {\rm OR}$). Let $m(\alpha)$ be the
unique $m < \omega$ such that 
$\rho_{m+1}({\cal M}^{\cal T}_{\alpha_{n(\alpha)}^*} || \eta(\alpha))
\leq \aleph_\alpha^H < 
\rho_m({\cal M}^{\cal T}_{\alpha_{n(\alpha)}^*} || \eta(\alpha))$. The following is
easy to verify.

\bigskip
{\bf Claim 2.} We may partition $C$ into finitely many sets $C_0$, ..., $C_k$, $0 \leq
k < \omega$, such that ${\rm sup}(C_l) \leq {\rm min}(C_{l+1})$ whenever $l<k$ and
such that for all $l \leq k$, if $\{ \alpha , \alpha' \} \subset C_l$ then
$n(\alpha) = n(\alpha')$, $\eta(\alpha) = \eta(\alpha')$, and $m(\alpha) =
m(\alpha')$.

\bigskip
In order to finish the proof of the Main Claim 
it therefore now suffices to find, for an arbitrary $l \leq k$, 
some $g \in K$, $g \colon \lambda
\rightarrow \lambda$ such that $g(\aleph_{\alpha+1}) > f(\aleph_{\alpha+1})$ for all
$\alpha \in C_l$.

Let us fix $l \leq k$. Let us write $n$, $\eta$, and $m$ for $n(\alpha)$,
$\eta(\alpha)$, and $m(\alpha)$, where $\alpha$ is any member of $C_l$.
Let $${\cal M} = {\cal M}_{\alpha^*_n}^{\cal T} || \eta.$$ 
Let $\lambda_l = {\rm sup}(\{ \aleph_\alpha^H \colon \alpha \in C_l \})$. We may
define $${\tilde \pi} \colon {\cal M} \rightarrow {\tilde {\cal M}} =
Ult_m({\cal M};\pi \upharpoonright \lambda_l).$$
By the argument of \cite{covering}, we shall have that ${\tilde {\cal M}}
\triangleleft K$. In particular, ${\tilde {\cal M}} \in K$.

Let us write $\gamma^-$ for the cardinal predecessor of $\gamma$ if $\gamma$ is a
successor cardinal (otherwise we let $\gamma^- = \gamma$).
Let ${\tilde \lambda}_l = {\rm sup}(\pi {\rm " } \lambda_l)$.
Let us define $g \colon {\tilde \lambda}_l \rightarrow {\tilde \lambda}_l$
as follows. Let $\gamma < {\tilde \lambda}_l$, and let us write $\gamma^-$
for $(\gamma^-)^{{\tilde {\cal M}}}$, i.e., $\gamma^-$ in the sense of ${\tilde {\cal
M}}$. 
Let $$\sigma_\gamma \colon {\tilde {\cal M}}_\gamma \cong h^{{\tilde {\cal
M}}}_m(\gamma^- \cup \{ p_{{\tilde {\cal M}}} \}) \prec_m {\tilde {\cal M}} {\rm , }$$
where 
${\tilde {\cal M}}_\gamma$ is transitive. We let $$g(\lambda) = 
(\gamma^-)^{+{\tilde {\cal M}}_\gamma}.$$
Notice that $g \in K$. We are left with having to verify that 
$g(\aleph_{\alpha+1}) > f(\aleph_{\alpha+1})$ for all
$\alpha \in C_l$.

Fix $\alpha \in C_l$. 
Let us assume without loss of generality that $\beta(\alpha) < \alpha_n^*$, the other
case being easier. Then $\eta = {\cal M}^{\cal T}_{\alpha_n^*} \cap {\rm OR}$, i.e.,
${\cal M} = {\cal M}^{\cal T}_{\alpha_n^*}$.
Consider $$\pi^{\cal T}_{\beta(\alpha) \alpha_n^*} \colon {\cal M}^{\cal
T}_{\beta(\alpha)} \rightarrow {\cal M}^{\cal T}_{\alpha_n^*}.$$
It is easy to verify that
$${\tilde {\cal M}}_{\aleph_{\alpha+1}} = Ult_m({\cal M}^{\cal
T}_{\beta(\alpha)};\pi \upharpoonright \pi^{-1}(\aleph_{\alpha})) {\rm , }$$
and that there is a map $$\varphi_\alpha \colon {\tilde {\cal M}}_{\aleph_{\alpha+1}}
\rightarrow {\tilde {\cal M}}$$ which is defined by 
$$\tau^{{\tilde {\cal M}}_{\aleph_{\alpha+1}}}({\vec \xi},p_{{\tilde 
{\cal M}}_{\aleph_{\alpha+1}}}) \mapsto 
\tau^{{\tilde {\cal M}}}({\vec \xi},p_{{\tilde 
{\cal M}}}) {\rm , }$$ where ${\vec \xi} < \aleph_\alpha$
and $\tau$ is an appropriate term, and which is just the
inverse of the collapsing map obtained from taking $h_m^{\tilde {\cal M}}(\aleph_\alpha
\cup \{ p_{\tilde {\cal M}} \})$. I.e., $\varphi_\alpha = \sigma_{\aleph_{\alpha+1}}$.

We now use the fact that $\alpha \in C$, i.e., that if $\beta+1 \in (0,\infty]_T$ then
we do not have that ${\rm crit}(E^{\cal T}_\beta) < \aleph_\alpha^H \leq \nu(E^{\cal
T}_\beta)$. This implies that $\pi^{\cal T}_{\beta(\alpha) \alpha_n^*}$
is the
inverse of the collapsing map obtained from taking $h_m^{\cal M}(\aleph_\alpha^H
\cup \{ p_{\cal M} \})$, and that $\aleph_{\alpha+1}^H \subset
{\cal M}^{\cal
T}_{\beta(\alpha)}$. In fact, $\aleph_{\alpha+1}^H =
\aleph_\alpha^{+{\bar K}}$, by \cite{covering}, $= \aleph_\alpha^{+{\cal M}^{\cal
T}_{\beta(\alpha)}}$. We therefore have that 
$$f(\aleph_{\alpha+1}) < {\rm sup}(\pi {\rm " } \aleph_{\alpha+1}^H) = 
\aleph_\alpha^{+{\tilde {\cal M}}_{\aleph_{\alpha+1}}}.$$
We have shown that $g(\aleph_{\alpha+1}) > f(\aleph_{\alpha+1})$.
\hfill $\square$

\bigskip
We would like to prove that if $2^{\aleph_1} < \aleph_\omega$, 
but $\aleph_\omega^{\aleph_0} >
\aleph_{\omega_1}$, then Projective Determinacy (or even $AD^{L({\mathbb R})}$) holds,
using Woodin's core model induction. This, however, would require a solution of
problem $\# 5$ of the list \cite{list}.

\end{document}